\documentclass[11pt,final]{article}%
\usepackage{stix2}
\usepackage{ifdraft}
\usepackage{amsfonts}
\usepackage{amsmath}
\usepackage{amsxtra}
\usepackage[amsmath,hyperref,thmmarks]{ntheorem}
\usepackage{makeidx}
\usepackage{graphicx}
\usepackage[colorlinks=true,bookmarksnumbered=true,pdfpagemode=None,final]%
{hyperref}
\usepackage{amssymb}%
\providecommand{\U}[1]{\protect\rule{.1in}{.1in}}
\ifdraft{}{\usepackage{version}}
\ifdraft{}{\excludeversion{nt}}
\def\mapsto{\DOTSB\mathchar"39AD }

\newcommand{\df}{\smash{\lower.12em\hbox{\textup{\tiny def}}}}

\usepackage[overload]{textcase}

\usepackage{xcolor}
\definecolor{bblue}{rgb}{0.0, 0.0, 0.6}

\usepackage{array}

\usepackage{microtype}

\usepackage[only,mapsfrom]{stmaryrd}

\usepackage{tikz}
\usepackage{tikz-cd}
\usetikzlibrary{matrix,arrows,positioning,decorations.pathmorphing}
\usetikzlibrary{decorations.markings,shapes.geometric}
\tikzset{commutative diagrams/column sep/Huge/.initial=24ex}
\tikzcdset{arrow style=math font}

\usepackage{wrapfig}

\usepackage{enumitem}
\setitemize[1]{label=$\diamond$\hspace{0.07in}}
\setlist{topsep=0.2em,itemsep=0.2em,parsep=0.2em}
\setdescription{font=\normalfont}

\usepackage{url}
\usepackage{verbatim}

\usepackage[notcite,notref]{showkeys}

\usepackage{mathtools}

\usepackage{titletoc,titlesec}
\contentsmargin{2.0em}
\dottedcontents{section}[2.3em]{}{1.8em}{1pc}
\usepackage[nottoc]{tocbibind}

\titleformat*{\section}{\LARGE\bfseries}
\titleformat*{\subsection}{\Large\itshape}
\titleformat*{\subsubsection}{\scshape}%\large
\titleformat*{\paragraph}{\itshape}

\usepackage{natbib}
\let\cite\citealt
\hypersetup{linkcolor=bblue,anchorcolor=bblue,citecolor=bblue,filecolor=bblue,urlcolor=bblue}
\hypersetup{pdftitle={Classification of Mumford--Tate groups},pdfauthor={J.S. Milne}}
\hypersetup{pdfsubject={Hodge theory},pdfkeywords={Hodge theory, Mumford--Tate groups}}

\usepackage[papersize={6.5in,10.0in},margin=0.5in,pdftex]{geometry}

\newcommand{\babstract}{\begin{abstract}}\newcommand\eabstract{\end{abstract}}
\newcommand{\bcomment}{\begin{comment}}\newcommand\ecomment{\end{comment}}
\newcommand{\bfootnotesize}{\begin{footnotesize}}\newcommand\efootnotesize{\end{footnotesize}}
\newcommand{\bquote}{\begin{quote}}\newcommand\equote{\end{quote}}
\newcommand{\bsmall}{\begin{small}}\newcommand\esmall{\end{small}}
\newcommand{\btable}{\begin{table}}\newcommand{\etable}{\end{table}}

\newcommand{\edocument}{\end{document}}

\newcommand{\tableoc}{\tableofcontents}%* removes "contents" from TOC
\newcommand{\textnf}{\textnormal}

\renewcommand{\Gamma}{\varGamma}
\renewcommand{\Pi}{\varPi}
\renewcommand{\Sigma}{\varSigma}
\def\1{{1\mkern-7mu1}}

\DeclareMathOperator{\ad}{ad}

\DeclareMathOperator{\End}{End}

\DeclareMathOperator{\Gal}{Gal}
\DeclareMathOperator{\GL}{GL}
\DeclareMathOperator{\Hom}{Hom}

\DeclareMathOperator{\MT}{MT}

\newcommand{\Hdg}{\mathsf{Hdg}}%Hodge structures

\theoremnumbering{arabic}
\theoremheaderfont{\scshape}
\RequirePackage{latexsym}
\theorembodyfont{\slshape}
\theoremseparator{.}
\newtheorem{X}{X}[section]

\newtheorem{corollary}[X]{Corollary}

\newtheorem{lemma}[X]{Lemma}

\newtheorem{proposition}[X]{Proposition}
\newtheorem{theorem}[X]{Theorem}

\theorembodyfont{\upshape}

\newtheorem{plain}[X]{}

\newtheorem{remark}[X]{Remark}

\theorembodyfont{\small}

\newtheorem*{note}{Notes}
\ifdraft{}{}
\theorembodyfont{\normalsize}
\theoremstyle{nonumberplain}
\theoremsymbol{\ensuremath{\color{lightgray}\blacksquare}}
\newtheorem{proof}{Proof}
\qedsymbol{\ensuremath{\color{lightgray}\blacksquare}}
\makeindex
\begin{document}

\title{Classification of the Mumford--Tate Groups of Rational Polarizable Hodge Structures}
\author{J.S. Milne}
\date{\today, v2.0}
\maketitle

\begin{abstract}
Let $G$ be the pro-algebraic group attached to the tannakian category of
polarizable rational Hodge structures. We show that the quotient of $G$ by its
derived group is the Serre group, the derived group of $G$ is the simply
connected covering of the adjoint group of $G$, and that the adjoint group $G$
is a product of specific simple algebraic groups. As the Mumford--Tate groups
are exactly the algebraic quotients of $G$, this also describes them.

\end{abstract}

\tableoc

\bigskip\label{m01} Mumford and Tate originally defined their algebraic groups
for complex abelian varieties. However, the group depends only on the Hodge
structure attached to the abelian variety, and the notion was soon extended to
all rational Hodge structures. The groups are of most interest when the Hodge
structure is polarizable.

\begin{plain}
\label{m02} The Mumford--Tate group of a rational Hodge structure is an
algebraic group $G$ over $\mathbb{Q}{}$ equipped with a cocharacter $\mu
\colon\mathbb{G}_{m}\rightarrow G_{\mathbb{C}{}}$. The weight $w(\mu)$ of
$\mu$ is the cocharacter $-\mu-\bar{\mu}$ of $G_{\mathbb{C}}$. In \S 3, we
obtain the following criterion: a pair $(G,\mu)$ is the Mumford--Tate group of
a polarizable rational Hodge structure if and only if it satisfies the
following conditions:

\begin{description}
\item[\textmd{mt1:}] the weight $w(\mu)$ of $\mu$ is defined over
$\mathbb{Q}{}$ and is central;

\item[\textmd{mt2:}] $\ad\mu(-1)$ is a Cartan involution of $(G/w(\mathbb{G}%
_{m}))_{\mathbb{R}}$;

\item[\textmd{mt3:}] $\mu$ generates $G$ (i.e., if $H\subset G$ is such that
$\mu(\mathbb{C}^{\times})\subset H(\mathbb{C})$, then $H=G$).
\end{description}
\end{plain}

\begin{plain}
\label{m03} The polarizable rational Hodge structures form a tannakian
category $\Hdg_{\mathbb{Q}{}}$ over $\mathbb{Q}$ with a canonical (forgetful)
fibre functor. The corresponding Tannaka group $G_{\mathrm{Hg}}$ is the
pro-algebraic group over $\mathbb{Q}$ having the Mumford--Tate groups as its
algebraic quotients. Thus understanding the Mumford--Tate groups amounts to
understanding $G_{\mathrm{Hg}}$. We obtain the following results:

\begin{itemize}
\item the quotient of $G_{\mathrm{Hg}}$ by its derived group is the
(well-known) Serre protorus $S$;

\item the derived group of $G_{\mathrm{Hg}}$ is simply connected, and hence is
the simply connected covering of the adjoint group of $G_{\mathrm{Hg}}$;

\item the simple factors of the adjoint group of $G_{\mathrm{Hg}}$ are the
groups of the form $(G)_{F/\mathbb{Q}}$, where $F$ is a totally real number
field and $G$ is a geometrically simple algebraic group over $F$ such that
$(G)_{F/\mathbb{Q}}(\mathbb{R})$ has a compact maximal torus.
\end{itemize}
\end{plain}

\noindent The article is largely expository because all of the intermediate
results have long been available in the literature.

\subsection{Notation and terminology}

All vector spaces are finite dimensional. Complex conjugation on $\mathbb{C}$
is denoted by $z\mapsto\bar{z}$ or $z\mapsto\iota z$. The terminology
concerning algebraic groups is that of \cite{milne2017}. In particular,
semisimple and reductive algebraic groups are connected, and an adjoint
algebraic group is a semisimple group with trivial centre. The centre of $G$
is denoted by $Z(G)$. When $K/k$ is a finite extension of fields and $G$ is an
algebraic group over $K$, we let $(G)_{K/k}$ denote the algebraic group over
$k$ obtained from $G$ by restriction of scalars.

\section{Definitions}

The \emph{Deligne torus} $\mathbb{S}{}$ is defined to be $(\mathbb{G}%
_{m})_{\mathbb{C}{}/\mathbb{R}{}}$. Thus%
\[
\mathbb{S}(\mathbb{R})=\mathbb{C}^{\times},\quad\mathbb{S}_{\mathbb{C}{}%
}\simeq\mathbb{G}_{m}\times\mathbb{G}_{m}.
\]
The map $\mathbb{S}{}(\mathbb{R}{})\rightarrow\mathbb{S}{}(\mathbb{C}{})$
induced by $\mathbb{R}{}\rightarrow\mathbb{C}{}$ is $z\mapsto(z,\bar{z})$.
There are homomorphisms%
\[
\begin{tikzcd}[column sep=large,row sep=0ex]
\mathbb{G}_{m}\arrow[r,"w"]
&\mathbb{S}\arrow[r,"t"]
&\mathbb{G}_{m},
&[-2.5em]t\circ w=-2,\\
\mathbb{R}^{\times}\arrow[r,"a\mapsto a^{-1}"]
&\mathbb{C}^{\times}\arrow[r,"z\mapsto{z\bar{z}}"]
&\mathbb{R}^{\times}.
\end{tikzcd}
\]

We denote the kernel of $t$ by $\mathbb{S}^{1}$. Thus $\mathbb{S}^{1}$ is a
one-dimensional torus over $\mathbb{R}$ with
\[
\mathbb{S}{}^{1}(\mathbb{R}{}) =\{z\in\mathbb{C}{}^{\times}\mid z\bar{z}=1\}
=\text{circle group $S^{1}$}.
\]
There is a canonical isomorphism
\[
\mathbb{S}{}/w(\mathbb{G}_{m})\rightarrow\mathbb{S}{}^{1},\quad
z\,\,(\mathrm{mod\,\mathbb{R}}^{\times})\mapsto z/\bar{z},
\]
with inverse $u\mapsto\sqrt{u}\,\,(\mathrm{mod\,\mathbb{R}{}}^{\times})$.

A homomorphism $h\colon\mathbb{S}\rightarrow G$ of real algebraic groups gives
rise to cocharacters%
\[%
\begin{array}
[c]{lll}%
\mu_{h}\colon\mathbb{G}_{m}\rightarrow G_{\mathbb{C}}, & z\mapsto
h_{\mathbb{C}{}}(z,1), & z\in\mathbb{G}_{m}(\mathbb{C})=\mathbb{C}^{\times},\\
w_{h}\colon\mathbb{G}_{m}\rightarrow G, & w_{h}=h\circ w & \text{(weight
homomorphism).}%
\end{array}
\]
The following formulas are useful,%
\begin{align*}
h_{\mathbb{C}{}}(z_{1},z_{2})  &  =\mu_{h}(z_{1})\cdot\bar{\mu}_{h}%
(z_{2});\quad h(z)=\mu_{h}(z)\cdot\overline{\mu_{h}(z)}\\
h(i)  &  =\mu_{h}(-1)\cdot w_{h}(i);\quad w_{h}=w(\mu_{h}).
\end{align*}

A \emph{Hodge structure} on a real vector space $V$ is a homomorphism
$h\colon\mathbb{S}\rightarrow\GL_{V}$. Such a homomorphism determines a
decomposition $V\otimes\mathbb{C}=\bigoplus V^{p,q}$, where $V^{p,q}$ is the
subspace on which $h(z)$ acts as $z^{-p}\cdot\iota z^{-q}$. A \emph{rational
Hodge structure} $(V,h)$ is a $\mathbb{Q}$-vector space $V$ together with a
Hodge structure $h$ on $V\otimes\mathbb{R}$ such that $w_{h}$ is defined over
$\mathbb{Q}{}$. The \emph{Tate Hodge structure} $\mathbb{Q}(m)$ is the
$\mathbb{Q}$-subspace $(2\pi i)^{m}\mathbb{Q}$ of $\mathbb{C}{}$ with $h(z)$
acting as multiplication by $(z\bar{z})^{m}$.

A \emph{polarization} of a rational Hodge structure $(V,h)$ of weight $m$ is a
morphism of Hodge structures
\[
\psi\colon V\otimes V\rightarrow\mathbb{\mathbb{Q}}(-m),\quad m\in\mathbb{Z},
\]
such that
\[
(x,y)\mapsto(2\pi i)^{m}\psi_{\mathbb{R}{}}(x,Cy)\colon V_{\mathbb{R}{}}\times
V_{\mathbb{R}{}}\rightarrow\mathbb{R}%
\]
is symmetric and positive definite. Here $C\overset{\df}{=}h(i)$ is the Weil operator.

\begin{note}
The conventions are those of \cite{deligne1979}.
\end{note}

\section{Cartan involutions and polarizations}

Let $G$ be a connected algebraic group over $\mathbb{R}$, and let $\sigma
_{0}\colon g\mapsto\bar{g}$ denote complex conjugation on $G(\mathbb{C})$ with
respect to $G$. A \emph{Cartan involution }%
\index{Cartan involution}
of $G$ is an involution $\theta$ of $G$ (as an algebraic group over
$\mathbb{R}{}$) such that the group%
\index{Gt@$G^{(\theta)}$}%
\[
G^{(\theta)}(\mathbb{R}{})=\{g\in G(\mathbb{C}{})\mid g=\theta(\bar{g})\}
\]
is compact. Then $G^{(\theta)}$ is a compact real form of $G_{\mathbb{C}}$,
and $\theta$ acts on $G(\mathbb{C})$ as $\sigma_{0}\sigma=\sigma\sigma_{0}$,
where $\sigma$ is complex conjugation on $G(\mathbb{C})$ with respect to
$G^{(\theta)}$.

A connected algebraic group $G$ over $\mathbb{R}{}$ has a Cartan involution if
and only if it has a compact real form, which is the case if and only if $G$
is reductive. Any two Cartan involutions of $G$ are conjugate by an element of
$G(\mathbb{R})$.

Let $C$ be an element of $G(\mathbb{R}{})$ whose square is central, so
$\ad(C)\overset{\df}{=}(g\mapsto CgC^{-1})$ is an involution. A $C$%
-\emph{polarization}%
\index{C-polarization@$C$-polarization}
on a real representation $V$ of $G$ is a $G$-invariant bilinear form
$\varphi\colon V\times V\rightarrow\mathbb{R}{}$ such that the form
$\varphi_{C}\colon(x,y)\mapsto\varphi(x,Cy)$ is symmetric and positive definite.

\begin{theorem}
\label{m1}If $\ad(C)$ is a Cartan involution of $G$, then every finite
dimensional real representation of $G$ carries a $C$-polarization; conversely,
if one \textup{faithful} finite dimensional real representation of $G$ carries
a $C$-polarization, then $\ad(C)$ is a Cartan involution.
\end{theorem}

\begin{proof}
An $\mathbb{R}{}$-bilinear form $\varphi$ on a real vector space $V$ defines a
sesquilinear form $\varphi^{\prime}\colon(u,v)\mapsto\varphi_{\mathbb{C}{}%
}(u,\bar{v})$ on $V(\mathbb{C}{})$, and $\varphi^{\prime}$ is hermitian (and
positive definite) if and only if $\varphi$ is symmetric (and positive definite).

Let $G\rightarrow\GL_{V}$ be a representation of $G$. If $\ad(C)$ is a Cartan
involution of $G$, then $G^{(\ad C)}(\mathbb{R}{})$ is compact, and so there
exists a $G^{(\ad C)}$-invariant positive definite symmetric bilinear form
$\varphi$ on $V$. Then $\varphi_{\mathbb{C}{}}$ is $G(\mathbb{C}{}%
)$-invariant, and so%
\[
\varphi^{\prime}(gu,(\sigma g)v)=\varphi^{\prime}(u,v),\quad\text{for all
}g\in G(\mathbb{C}{})\text{, }u,v\in V_{\mathbb{C}{}},
\]
where $\sigma$ is the complex conjugation on $G_{\mathbb{C}{}}$ with respect
to $G^{(\ad C)}$. Now $\sigma g=\ad(C)(\bar{g})=\ad(C^{-1})(\bar{g})$, and so,
on replacing $v$ with $C^{-1}v$ in the equality, we find that%
\[
\varphi^{\prime}(gu,(C^{-1}\bar{g}C)C^{-1}v)=\varphi^{\prime}(u,C^{-1}%
v),\quad\text{for all }g\in G(\mathbb{C}{})\text{, }u,v\in V_{\mathbb{C}{}}.
\]
In particular, $\varphi(gu,C^{-1}gv)=\varphi(u,C^{-1}v)$ when $g\in
G(\mathbb{R}{})$ and $u,v\in V$. Therefore, $\varphi_{C^{-1}}$ is
$G$-invariant. As $(\varphi_{C^{-1}})_{C}=\varphi$, we see that $\varphi$ is a
$C$-polarization.

For the converse, one shows that, if $\varphi$ is a $C$-polarization on a
faithful representation, then $\varphi_{C}$ is invariant under $G^{(\ad
C)}(\mathbb{R}{})$, which is therefore compact.
\end{proof}

\begin{plain}
\label{m2}Let $G$ be an algebraic group over $\mathbb{Q}{}$, and let $C$ be an
element of $G(\mathbb{R}{})$ whose square is central. A $C$%
-\emph{polarization} on a $\mathbb{Q}{}$-representation $V$ of $G$ is a
$G$-invariant bilinear form $\varphi\colon V\times V\rightarrow\mathbb{Q}{}$
such that $\varphi_{\mathbb{R}{}}$ is a $C$-polarization on $V_{\mathbb{R}{}}%
$. In order to show that a $\mathbb{Q}{}$-representation $V$ of $G$ is
polarizable, it suffices to check that $V_{\mathbb{R}{}}$ is polarizable. We
prove this when $C^{2}$ acts as $+1$ or $-1$ on $V$, which are the only cases
we shall need. Let $P(\mathbb{Q})$ (resp.\ $P(\mathbb{R}{})$) denote the space
of $G$-invariant bilinear forms on $V$ (resp.\ on $V_{\mathbb{R}{}}$) that are
symmetric when $C^{2}$ acts as $+1$ or skew-symmetric when it acts as $-1$.
Then $P(\mathbb{R}{})=P(\mathbb{Q}\otimes_{\mathbb{Q}}\mathbb{R})$. The
$C$-polarizations of $V_{\mathbb{R}{}}$ form an open subset of $P(\mathbb{R}%
{})$, whose intersection with $P(\mathbb{Q}{})$ consists of the $C$%
-polarizations of $V$.
\end{plain}

\begin{note}
Theorem 2.1 is \cite{deligne1972}, 2.8. The exposition follows
\cite{milne2005}, 1.20.
\end{note}

\section{Mumford--Tate groups}

Let $(V,h)$ be a rational Hodge structure. Following \cite{deligne1972}, 7.1,
we define the \emph{Mumford--Tate group}%
\index{Mumford--Tate group}
of $(V,h)$ to be the smallest algebraic subgroup $G$ of $\GL_{V}$ such that
$G_{\mathbb{R}{}}\supset h(\mathbb{S}{})$. We usually regard the Mumford--Tate
group as a pair $(G,h)$. Note that $G$ is connected, because otherwise we
could replace it with its neutral component.

The rational Hodge structures form a tannakian category over $\mathbb{Q}$. Let
$(V,h)$ be a rational Hodge structure, and let $\langle V,h\rangle^{\otimes}$
be the tannakian subcategory generated by $(V,h)$. The Mumford--Tate group of
$(V,h)$ is the algebraic group attached to $\langle V,h\rangle^{\otimes}$ and
the forgetful fibre functor.

The \emph{special Mumford--Tate group} of $(V,h)$ is defined to be the
smallest algebraic subgroup $G^{1}$ of $\GL_{V}$ such that $G_{\mathbb{R}{}%
}^{1}\supset h(\mathbb{S}{}^{1})$. It is a subgroup of the Mumford-Tate group
$G$, and $G=G^{1}\cdot w_{h}(\mathbb{G}_{m})$.

Let $G$ be a connected algebraic group over $\mathbb{Q}$ and $h$ a
homomorphism $\mathbb{S}{}\rightarrow G_{\mathbb{R}{}}$. Consider the
following conditions\footnote{These are the conditions SV4 and SV2* of
\cite{milne2005}, which, for a reductive group, coincide with the conditions
(2.1.1.4) and (2.1.1.5) of \cite{deligne1979}.} on $h$:

\begin{description}
\item[\textmd{MT1:}] the map $w_{h}\colon\mathbb{G}_{m\mathbb{R}{}}\rightarrow
G_{\mathbb{R}{}}$ is defined over $\mathbb{Q}$ and $w_{h}(\mathbb{G}%
_{m})\subset Z(G)$;

\item[\textmd{MT2:}] $\ad(h(i))$ is a Cartan involution of $(G/w_{h}%
(\mathbb{G}_{m}))_{\mathbb{R}}$.
\end{description}

\noindent Note that (MT1) implies that $G/w_{h}(\mathbb{G}_{m})$ is an
algebraic group over $\mathbb{Q}$ and that (MT2) implies that $G$ is reductive.

\begin{theorem}
\label{m3a}A pair $(G,h)$ as above is the Mumford--Tate group of a polarizable
rational Hodge structure if and only if it satisfies (MT1,2) and $h$ generates
$G$.
\end{theorem}

This combines the next two propositions.

\begin{proposition}
\label{m3}A pair $(G,h)$ as above is the Mumford--Tate group of a rational
Hodge structure if and only if $h$ satisfies (MT1) and $h$ generates $G$.
\end{proposition}

\begin{proof}
If $(G,h)$ is the Mumford--Tate group of a Hodge structure $(V,h)$, then
certainly $h$ generates $G$. The weight homomorphism $w_{h}$ is defined over
$\mathbb{Q}{}$ because $(V,h)$ is a rational Hodge structure. Let $Z(w_{h})$
denote the centralizer of $w_{h}$ in $G$. For any $a\in\mathbb{R}^{\times}$,
$w_{h}(a)\colon V_{\mathbb{R}{}}\rightarrow V_{\mathbb{R}{}}$ is a morphism of
real Hodge structures, and so it commutes with the action of $h(\mathbb{S}{}%
)$. Hence $h(\mathbb{S}{})\subset Z(w_{h})_{\mathbb{R}{}}$. As $h$ generates
$G$, this implies that $Z(w_{h})=G$.

Conversely, suppose that $(G,h)$ satisfies the conditions. For any faithful
representation $\rho\colon G\rightarrow\GL_{V}$ of $G$, the pair
$(V,h\circ\rho)$ is a rational Hodge structure, and $(G,h)$ is its
Mumford--Tate group.
\end{proof}

\begin{proposition}
\label{m4}Let $(G,h)$ be the Mumford--Tate group of a rational Hodge structure
$(V,h)$. Then $(V,h)$ is polarizable if and only if $h$ satisfies (MT2).
\end{proposition}

\begin{proof}
Let $C=h(i)$. For notational convenience, assume that $(V,h)$ has a single
weight $m$. Let $G^{1}$ be the special Mumford--Tate group of $(V,h)$. Then
$C\in G^{1}(\mathbb{R}{})$, and a pairing $\psi\colon V\times V\rightarrow
\mathbb{Q}(-m)$ is a polarization of the Hodge structure $(V,h)$ if and only
if $(2\pi i)^{m}\psi$ is a $C$-polarization of $V$ for $G^{1}$ in the sense of
\S 2. It follows from (\ref{m1}) and (\ref{m2}) that a polarization $\psi$ for
$(V,h)$ exists if and only if $\ad(C)$ is a Cartan involution of
$G_{\mathbb{R}{}}^{1}$. Now $G^{1}\subset G$ and the quotient map
$G^{1}\rightarrow G/w_{h}(\mathbb{G}_{m})$ is an isogeny, and so $\ad(C)$ is a
Cartan involution of $G^{1}$ if and only if it is a Cartan involution of
$G/w_{h}(\mathbb{G}_{m})$.
\end{proof}

\begin{corollary}
\label{m5}The Mumford--Tate group of a polarizable rational Hodge structure is reductive.
\end{corollary}

\begin{proof}
Immediate consequence of Proposition \ref{m4}.
\end{proof}

There is a canonical homomorphism $h_{\mathrm{Hg}}\colon\mathbb{S}%
{}\rightarrow(G_{\mathrm{Hg}})_{\mathbb{R}{}}$ corresponding to the functor
$-\otimes\mathbb{R}{}$ from polarizable rational Hodge structures to
polarizable real Hodge structures$.$

\begin{corollary}
\label{m20}For any reductive group $G$ over $\mathbb{Q}{}$ and homomorphism
$h\colon\mathbb{S}{}\rightarrow G_{\mathbb{R}}$ satisfying (MT2,4), there is a
unique homomorphism $\rho\colon G_{\mathrm{Hg}}\rightarrow G$ such that
$\rho_{\mathbb{R}{}}\circ h_{\mathrm{Hg}}=h$.
\end{corollary}

\begin{proof}
Immediate consequence of Theorem \ref{m3a}.
\end{proof}

\begin{remark}
\label{m24}Let $(V,h)$ be a rational Hodge structure, and let $\mu=\mu_{h}$.
Then $h(z)=\mu_{h}(z)\cdot\overline{\mu_{h}(z)}$ and so $\mu_{h}$ determines
$h$. A cocharacter $\mu$ of $G_{\mathbb{C}{}}$ is of the form $\mu_{h}$ if and
only if $\mu$ commutes with $\bar{\mu}$. The Mumford--Tate group of $(V,h)$ is
the smallest algebraic subgroup $G$ of $\GL_{V}$ such that $G_{\mathbb{C}{}%
}\supset\mu_{h}(\mathbb{G}_{m})$. As $h(i)=\mu(-1)\cdot w_{h}(i)$, we see that
(\ref{m02}) is simply a restatement of Theorem \ref{m3a}.
\end{remark}

From now on, we say \textquotedblleft$(G,h)$ or $(G,\mu)$ is a Mumford--Tate
group\textquotedblright\ to mean that the pair is the Mumford--Tate group of a
\textit{polarizable} rational Hodge structure. We say that $G$ is a
Mumford-Tate group if there exists an $h$ such that $(G,h)$ is a Mumford--Tate group.

\begin{note}
Theorem 3.1 is Proposition 1.6 of \cite{milne1994}. The exposition follows
\cite{milne2013}, \S 6.
\end{note}

\section{Tori as Mumford--Tate groups}

A number field $E$ is a \emph{CM field} if it is a totally imaginary quadratic
extension of a totally real field. Let $\mathbb{Q}{}^{\mathrm{cm}}$ be the
union of the CM-subfields of $\mathbb{Q}{}^{\mathrm{al}}$. Then $\mathbb{Q}%
^{\mathrm{cm}}$ is the largest Galois extension of $\mathbb{Q}{}$ in
$\mathbb{Q}{}^{\mathrm{al}}$ such that complex conjugation is in the centre of
$\Gal(\mathbb{Q}{}^{\mathrm{cm}}/\mathbb{Q}{})$.

\begin{lemma}
\label{m6a}Let $T$ be a torus over $\mathbb{Q}{}$ and $\mu$ a cocharacter of
$T$ over $\mathbb{Q}{}^{\mathrm{al}}$. The following conditions on $\mu$ are equivalent:

\begin{enumerate}
\item the weight of $\mu$ is defined over $\mathbb{Q}{}$ and $\mu$ is defined
over a CM field;

\item for all $\sigma\in\Gal(\mathbb{Q}^{\mathrm{al}}/\mathbb{Q})$,
\[
(\sigma-1)(\iota+1)\mu=0=(\iota+1)(\sigma-1)\mu.
\]

\end{enumerate}
\end{lemma}

\begin{proof}
The first equality in (b) says that the weight of $\mu$ is defined over
$\mathbb{Q}{}$, and then the second says that $\sigma\iota\mu=\iota\sigma\mu$
for all $\sigma\in\Gal(\mathbb{Q}^{\mathrm{al}}/\mathbb{Q})$, i.e., that $\mu$
is defined over $\mathbb{Q}{}^{\mathrm{cm}}$.
\end{proof}

The equivalent conditions of the lemma are called the \emph{Serre condition}.
When $T$ is split by a CM field, the Serre condition simply says that the
weight of $\mu$ is defined over $\mathbb{Q}{}$.

A rational Hodge structure $(V,h)$ is said to be of \emph{CM-type \/} if it is
polarizable and its Mumford--Tate group is commutative (hence a torus). When
$(V,h)$ is simple, this means that $\End(V,h)$ is either a CM-field or
$\mathbb{Q}$.

We have the following criterion.

\begin{proposition}
\label{m6}Let $T$ be a torus over $\mathbb{Q}{}$ and $\mu$ a cocharacter of
$T$ over $\mathbb{Q}{}^{\mathrm{al}}$. Then $(T,\mu)$ is a Mumford--Tate group
if and only if

\begin{enumerate}
\item the weight of $\mu$ is defined over $\mathbb{Q},$

\item $T$ is split by a CM field,

\item $\mu$ generates $T$.
\end{enumerate}
\end{proposition}

\begin{proof}
When $\mu$ generates $T$, the condition (a)+(b) is equivalent to (a) of Lemma
\ref{m6a}; on the other hand, the condition (mt1)+(mt2) is equivalent to (b)
of Lemma \ref{m6a}. Thus, the proposition is a restatement of Theorem
\ref{m3a} for the case of tori.
\end{proof}

Let $E$ be a CM subfield of $\mathbb{Q}^{\mathrm{al}}$. Then $(\mathbb{G}%
_{m})_{E/\mathbb{Q}}$ is a torus with character group $\mathbb{Z}%
^{\Hom(E,\mathbb{\mathbb{Q}}{}^{\mathrm{al}})}$, and we define $S^{E}$ to be
the quotient of $(\mathbb{G}_{m})_{E/\mathbb{Q}}$ such that
\[
X^{\ast}(S^{E})= \{\lambda\in\mathbb{Z}^{\Hom(E,\mathbb{Q}^{\mathrm{al}})}%
\mid\lambda(\sigma)+\lambda(\iota\circ\sigma)=\text{constant},\,\,\,\sigma
\in\Hom(E,\mathbb{Q}^{\mathrm{al}})\}.
\]
Define $\mu^{E}$ to be the cocharacter of $S^{E}$ such that
\[
\langle\lambda,\mu^{E}\rangle=\lambda(\sigma_{0}),\quad\text{all }\lambda\in
X^{\ast}(S^{E}),
\]
where $\sigma_{0}$ is the given embedding of $E$ into $\mathbb{Q}%
{}^{\mathrm{al}}$. If $E\subset E^{\prime}\subset\mathbb{Q}{}^{\mathrm{al}}$,
then the norm map defines a homomorphism $S^{E^{\prime}}\rightarrow S^{E}$
carrying $\mu^{E^{\prime}}$to $\mu^{E}$. We set
\[
(S,\mu_{\text{can}})=\varprojlim(S^{E},\mu^{E}).
\]
The pair $(S,\mu_{\text{can}})$ is called the \emph{Serre group}. Note that
$X^{\ast}(S)$ is the set of all locally constant functions $\lambda
\colon\Gal(\mathbb{Q}^{\mathrm{cm}}/\mathbb{Q})\rightarrow\mathbb{Z}$ such
that
\[
\lambda(\sigma)+\lambda(\iota\circ\sigma)=-m
\]
for some integer $m$ (called the \emph{weight\/} of $\lambda$).

\begin{theorem}
\label{m8} (a) The cocharacter $\mu^{E}$ of $S^{E}$ satisfies the Serre condition.

(b) Let $T$ be a torus over $\mathbb{Q}{}$ and $\mu$ a cocharacter satisfying
the Serre condition. Then there is a unique homomorphism $\rho_{\mu}\colon
S\rightarrow T$ such that $(\rho_{\mu})_{\mathbb{Q}}\circ\mu_{\text{can}}=\mu$.

(c) We have
\begin{equation}
(S,\mu_{\text{can}})=\varprojlim(T,\mu), \label{qm1}%
\end{equation}
where the limit runs over the pairs $(T,\mu)$ such that $\mu$ satisfies the
Serre condition and generates $T$.
\end{theorem}

\begin{proof}
(a) The weight of $\mu$ is defined over $\mathbb{Q}$ and $T$ is split by a CM field.

(b) For $\chi\in X^{\ast}(T)$ and $\sigma\in\Gal(\mathbb{Q}{}^{\mathrm{al}%
}/\mathbb{Q}{})$, define%
\[
f_{\chi}(\sigma)=\langle\sigma^{-1}\chi,\mu\rangle.
\]
Then%
\begin{equation}
\chi\mapsto f_{\chi}\colon X^{\ast}(T)\rightarrow X^{\ast}(S) \label{qm2}%
\end{equation}
is a $\Gal(\mathbb{Q}{}^{\mathrm{al}}/\mathbb{Q}{})$-equivariant homomorphism,
and so corresponds to a homomorphism%
\[
\rho\colon S\rightarrow T.
\]
For $\sigma_{0}$ the given inclusion of $E$ into $\mathbb{Q}{}^{\mathrm{al}}$,
$f_{\chi}(\sigma_{0})=\langle\mu,\chi\rangle$, i.e., $\langle\mu
_{\mathrm{can}},f_{\chi}\rangle=\langle\mu,\chi\rangle$, which shows that
$(\rho_{\mu})_{\mathbb{Q}}\circ\mu_{\text{can}}=\mu$.

(c) For every $(T,\mu)$, we have defined an injective homomorphism
(\ref{qm2}), and $X^{\ast}(S)$ is union of their images.
\end{proof}

\begin{corollary}
\label{m7} The polarizable rational Hodge structures of CM-type form a
tannakian category, and $S$ is the pro-algebraic group attached to the
forgetful fibre functor.
\end{corollary}

\begin{proof}
For any rational Hodge structures $X$ and $Y$, $\MT(X\oplus Y)\subset
\MT(X)\times\MT(Y)$, and so $X\oplus Y$ is of CM type if $X$ and $Y$ are. The
category of polarizable rational Hodge structures of CM-type is the directed
union of the categories $\langle X\rangle^{\otimes}$ with $X$ of CM-type, and
is therefore tannakian. Correspondingly, the pro-algebraic group attached to
the forgetful fibre functor is $\varprojlim(\MT(X),\mu_{X})$, which, according
to \ref{m6} and (\ref{qm1}), is equal to $(S,\mu_{\text{can}})$.
\end{proof}

The functor sending a rational Hodge structure $(V,h)$ to the real Hodge
structure $(V\otimes\mathbb{R}{},h)$ defines a homomorphism $h_{\mathrm{can}%
}\colon\mathbb{S}{}\rightarrow S$; its associated cocharacter is
$\mu_{\mathrm{can}}$.

\begin{note}
Everything in this section has been known to the experts since the 1960s ---
see, for example, \cite{serre1968}. For a detailed account, see my notes
\textit{Complex Multiplication}.
\end{note}

\section{Semisimple groups as Mumford-Tate groups}

Let $G$ be an algebraic group over $\mathbb{R}$ and $h\colon\mathbb{S}%
{}\rightarrow G$ a homomorphism of weight $0$. Then $\iota\mu_{h}=-\mu_{h}$,
and so $\mu_{h}\colon(\mathbb{G}_{m})_{\mathbb{C}{}}\rightarrow G_{\mathbb{C}%
{}}$ arises from a homomorphism $u\colon\mathbb{G}_{m}\rightarrow G$ over
$\mathbb{R}$. As $h$ has weight $0$, it factors through $\mathbb{S}%
{}/w(\mathbb{G}_{m})$, and $u$ is the composite%
\[
\mathbb{S}{}^{1}\simeq\mathbb{S}{}/w(\mathbb{G}_{m}%
)\overset{h}{\longrightarrow}G.
\]
In this way we get a one-to-one correspondence between homomorphisms
$h\colon\mathbb{S}{}\rightarrow G$ of weight $0$ and homomorphisms
$u\colon\mathbb{S}{}^{1}\rightarrow G$. If $h\leftrightarrow u,$ then
\begin{align*}
h(z)  &  =u(z/\bar{z}),\quad z\in\mathbb{C}^{\times},\\
u(z)  &  =h(\sqrt{z}),\quad z\in U^{1}.
\end{align*}
Note that $h(i)=u(-1)$, so (MT2) becomes the condition that $\ad u(-1)$ is a
Cartan involution.

\begin{lemma}
\label{m9a}Let $G$ be a simple algebraic group over $\mathbb{R}{}$ (so, in
particular, adjoint). If $G$ is an inner form of its compact form, then it is
geometrically simple.
\end{lemma}

\begin{proof}
If $G_{\mathbb{C}}$ is not simple, say, $G_{\mathbb{C}{}}=G_{1}\times G_{2}$,
then $G=(G_{1})_{\mathbb{C}{}/\mathbb{R}{}}$ and any inner form of $G$ is also
the restriction of scalars of a complex group, but such a group cannot be
compact (look at a subtorus).
\end{proof}

\begin{proposition}
\label{m9}Let $G$ be a simple algebraic group over $\mathbb{R}{}$. Then $G$
admits a homomorphism $u\colon\mathbb{S}{}^{1}\rightarrow G$ such that $\ad
u(-1)$ is a Cartan involution if and only if $G$ contains a compact maximal
torus (in which case, $G$ is geometrically simple).
\end{proposition}

\begin{proof}
$\Longrightarrow$: Any maximal torus containing $u(\mathbb{S}{}^{1})$ is compact.

$\Longleftarrow$: Let $T$ be a compact maximal torus of $G_{\mathbb{R}{}}$.
Choose a maximal compact subgroup of $G_{\mathbb{R}{}}$ containing $T$, and
let $\theta$ be the corresponding Cartan involution. A root of $(G,T)$ is
compact or noncompact according as $\theta$ acts as $-1$ or $+1$ on it. A
homomorphism $u$ such that $\langle\alpha,u\rangle$ is even or odd according
as $\alpha$ is compact or noncompact has the property that $\ad(u(-1))$ is a
Cartan involution of $G_{\mathbb{R}}$.
\end{proof}

\begin{remark}
\label{m10}It is possible to read off from the classification of geometrically
simple algebraic groups over $\mathbb{R}$, a list of the groups satisfying the
equivalent conditions of Proposition \ref{m9}.
\end{remark}

\begin{theorem}
\label{m10a}An adjoint algebraic group $G$ over $\mathbb{Q}{}$ is a
Mumford--Tate group if and only if $G_{\mathbb{R}{}}$ contains a compact
maximal torus.
\end{theorem}

\begin{proof}
As $G$ is adjoint, to give a homomorphism $h\colon\mathbb{S}{}\rightarrow
G_{\mathbb{R}{}}$ satisfying (MT1,2) is the same as giving a homomorphism
$u\colon\mathbb{S}{}^{1}\rightarrow G_{\mathbb{R}{}}$ such that $\ad u(-1)$ is
a Cartan involution. If $G_{\mathbb{R}{}}$ contains a maximal torus, then the
proof of \ref{m9} shows how to construct such a $u$. A general $u$ will
generate $G$, and so $G$ is a Mumford--Tate group by Theorem \ref{m3a}.
Conversely, if $G$ is a Mumford--Tate group, then Proposition \ref{m9} shows
that $G_{\mathbb{R}{}}$ contains a compact maximal torus.
\end{proof}

\begin{remark}
\label{m10b} Let $G$ be an adjoint algebraic group over $\mathbb{Q}$. Then $G$
has a canonical decomposition%
\[
G=(G_{1})_{F_{1}/\mathbb{Q}}\times\cdots\times(G_{n})_{F_{n}/\mathbb{Q}},
\]
where each $F_{i}$ is a subfield of $\mathbb{Q}{}^{\mathrm{al}}$ and ${}G_{i}$
is geometrically simple (\cite{milne2017}, 24.4). If the simple factors of
$(G_{i})_{F_{i}/\mathbb{Q}}$ over $\mathbb{R}{}$ are geometrically simple,
then $F_{i}$ is totally real. Thus, $G$ is a Mumford-Tate group if and only if
its simple factors (over $\mathbb{Q}$) are the groups of the form
$(H)_{F/\mathbb{Q}{}}$, where $F$ is a totally real number field and $H$ is a
geometrically simple group over $F$ such that, for every $\rho\colon
F\rightarrow\mathbb{R}$, $\rho H$ is on the list hinted at in \ref{m10}.
\end{remark}

\begin{note}
In their 2012 monograph, Green, Griffiths, and Kerr (IV.A.3) claim to show
that an adjoint group over $\mathbb{Q}{}$ is a Mumford--Tate group if and only
if it has an anisotropic maximal torus. Patrikis pointed out that this
statement is false and gave the correct statement.
\end{note}

\section{The classification}

We describe the structure of the pro-algebraic group $G_{\mathrm{Hg}}$
attached to the tannakian category of polarizable rational Hodge structures
and the forgetful fibre functor.

\begin{lemma}
[\cite{mumford1969}, \textnf{p.~348}]\label{m18}Let $G$ be a connected
algebraic group over $\mathbb{Q}$ and $T$ a maximal torus in $G_{\mathbb{R}}$.
Then there exists a maximal torus $T_{0}$ in $G$ and an $a\in G(\mathbb{R}{})$
such that $T_{0\mathbb{R}}=aTa^{-1}$.
\end{lemma}

\begin{proof}
According to the real approximation theorem (\cite{milne2017}, 25.70).We use
that $G(\mathbb{Q})$ is dense in $G(\mathbb{R})$ (real approximation theorem,
\cite{milne2017}, 25.70). If $a\in T(\mathbb{R}{})$ is a regular element, then
$T$ is the centralizer of $a$, and $a$ has an open neighbourhood $U$ in
$G(\mathbb{R}{})$ such that the centralizer of every $a^{\prime}\in U$ is a
conjugate of $T$. If $a^{\prime}\in U\cap G(\mathbb{Q}{})$, then the
centralizer of $a^{\prime}$ is a conjugate of $T$ defined over $\mathbb{Q}{}$,
as required.
\end{proof}

\begin{proposition}
[\cite{mumford1969}, \textnf{p.~348}]\label{m15} Let $G$ be an adjoint group
over $\mathbb{Q}{}$ and $h\colon\mathbb{S}{}\rightarrow G_{\mathbb{R}{}}$ a
homomorphism satisfying (MT1,2). There exists a $g\in G(\mathbb{R}{})$ and a
torus $T_{0}\subset G$ such that $\ad(g)\circ h$ factors through
$T_{0\mathbb{R}{}}$.
\end{proposition}

\begin{proof}
Let $K$ be the centralizer of $h$ in $G_{\mathbb{R}{}}$ (so $K$ is an
algebraic subgroup of $G_{\mathbb{R}{}}$). Let $T$ be a maximal torus of $K$.
As $h(\mathbb{\mathbb{S}{}})$ is contained in the centre of $K$,
$h(\mathbb{\mathbb{S}{}})\subset T$. If $T^{\prime}$ is a torus in
$G_{\mathbb{R}{}}$ containing $T$, then $T^{\prime}$ centralizes $h$ and so
$T^{\prime}\subset K$; therefore $T$ is maximal in $G_{\mathbb{R}{}}$.
According to the lemma, there exists a maximal torus $T_{0}$ of $G$such that
$T_{0\mathbb{R}{}}=gTg^{-1}$ for some $g\in G(\mathbb{R}{})$. Now $\ad(g)\circ
h$ factors through $T_{0\mathbb{R}{}}$.
\end{proof}

\begin{proposition}
\label{m12} Let $G$ be a reductive group over a field $k$ of characteristic
zero, and let $L$ be a finite Galois extension of $k$ splitting $G$. Let
$G^{\prime}\rightarrow G^{\mathrm{der}}$ be a finite covering of the derived
group of $G$. Then there exists a central extension
\[
1\rightarrow N\rightarrow G_{1}\rightarrow G\rightarrow1
\]
such that $G_{1}$ is reductive, $N$ is a product of copies of $(\mathbb{G}%
_{m})_{L/k}$, and
\[
(G_{1}^{\mathrm{der}}\rightarrow G^{\mathrm{der}})=(G^{\prime}\rightarrow
G^{\mathrm{der}}).
\]

\end{proposition}

\begin{proof}
See \cite{milneS1982}, 3.1.
\end{proof}

\begin{theorem}
[\cite{milne1994}, 1.28]\label{m16} Let $H$ be a semisimple algebraic group
over $\mathbb{Q}$ and $\bar{h}\colon\mathbb{S}/\mathbb{G}_{m}\rightarrow
H_{\mathbb{R}}^{\mathrm{ad}}$ a homomorphism such that $\ad(\bar{h}(i))$ is a
Cartan involution. Then there exists a reductive group $G$ with
$G^{\mathrm{der}}=H$ and a homomorphism $h\colon\mathbb{S}\rightarrow
G_{\mathbb{R}}$ lifting $\bar{h}$ and satisfying (MT1,2).
\end{theorem}

\begin{proof}
Assume first that $\bar{h}$ is \textquotedblleft special\textquotedblright,
i.e., that it factors through $T_{\mathbb{R}}$ for some maximal torus $T$ in
$H^{\mathrm{ad}}$. The hypothesis on $h$ implies that $T_{\mathbb{R}}$ is
anisotropic, and so $T$ splits over a CM-field $L$, which we may choose to be
Galois over $\mathbb{Q}$. According to Proposition \ref{m12}, there exists a
central extension defined over $\mathbb{Q}$
\[
1\longrightarrow N\longrightarrow G\longrightarrow H^{\mathrm{ad}%
}\longrightarrow1
\]
such that $G^{\mathrm{der}}=H$ and $N$ is a product of copies of
$(\mathbb{G}_{m})_{L/\mathbb{Q}}$. There is a maximal torus $T^{\prime}\subset
G$ mapping onto $T$ (\cite{milne2017}, 17.20). Since $T^{\prime}$ is its own
centralizer, it contains $N$, which is therefore the kernel of $T^{\prime
}\rightarrow T$. Hence $X_{\ast}(T^{\prime})\rightarrow X_{\ast}(T)$ is
surjective, and we can choose $\mu\in X_{\ast}(T^{\prime})$ mapping to
$\mu_{\bar{h}}\in X_{\ast}(T)$. The weight $w(h)\overset{\df}{=}-\mu-\iota\mu$
of $\mu$ lies in $X_{\ast}(N)$. Because $X_{\ast}(N)$ is an induced Galois
module, its cohomology groups are zero; in particular, the zeroth Tate group
\[
H_{\text{Tate}}^{0}(\Gal(\mathbb{C}/\mathbb{R}),X_{\ast}(N))\overset{\df}{=}%
\frac{X_{\ast}(N)^{\Gal(\mathbb{C}/\mathbb{R})}}{(\iota+1)X_{\ast}(N)}=0.
\]
Clearly $\iota w=w$, and so there exists a $\mu_{0}\in X_{\ast}(N)$ such that
$(\iota+1)\mu_{0}=w$. When we replace $\mu$ with $\mu+\mu_{0}$, then we find
that the weight becomes 0; in particular, it is defined over $\mathbb{Q}$.
Choose $h$ so that $h(z)=\mu(z)\cdot\overline{\mu(z)}$.

For a general $\bar{h}$, there will exist a $\bar{g}\in H^{\mathrm{ad}%
}(\mathbb{R})$ such that $\ad\bar{g}\circ\bar{h}$ is special (\ref{m15}).
Construct $G$ and $h$ as in the last paragraph corresponding to $\ad\bar
{g}\circ\bar{h}$. Because $H^{1}(\mathbb{R},N)=H^{1}(L\otimes_{\mathbb{Q}%
}\mathbb{R},\mathbb{G}_{m})=0$, the element $\bar{g}$ will lift to an element
$g\in G(\mathbb{R})$, and we take the pair $(G,\ad(g^{-1})\circ h)$.

For the pair $(G,h)$ we have constructed, the centre of $G$ is split by a
CM-field, $h$ satisfies (MT1), and $\ad(h(i))$ is a Cartan involution on
$G^{\mathrm{ad}}$. Let $T$ be the subtorus of $G/G^{\mathrm{der}}$ generated
by $h$. Then $T_{\mathbb{R}}$ is anisotropic, and when we replace $G$ with the
inverse image of $T$, we obtain a pair $(G,h)$ satisfying (MT1,2).
\end{proof}

\begin{theorem}
\label{m17}

\begin{enumerate}
\item The quotient of $G_{\mathrm{Hg}}$ by its derived group is the Serre group.

\item A semisimple algebraic group $G$ over $\mathbb{Q}$ is a quotient of
$G_{\mathrm{Hg}}^{\mathrm{der}}$ if and only if $G^{\mathrm{ad}}$ is a
Mumford--Tate group.

\item The adjoint group of $G_{\mathrm{Hg}}$ is a product of groups of the
form $(G)_{F/\mathbb{Q}{}}$ with $F$ a totally real number field and $G$ an
algebraic group over $F$ such that, for all embeddings $\rho$ of $F$ into
$\mathbb{R}{}$, $\rho G$ is a simple algebraic group over $\mathbb{R}{}$ with
a compact maximal torus.
\end{enumerate}
\end{theorem}

\begin{proof}
(a) A polarizable rational Hodge structure is of CM-type if and only if
$G_{\mathrm{Hg}}^{\mathrm{der}}$ acts trivially on it. Now \ref{m7} implies
that the inclusion of the category of CM Hodge structures into the full
category of polarizable rational Hodge structures induces an isomorphism
$G/G^{\mathrm{der}}\rightarrow S$.

(b) Let $G$ be a quotient of $G_{\mathrm{Hg}}^{\mathrm{der}}$. The image of
$Z(G_{\mathrm{Hg}}^{\mathrm{der}})$ in $G$ is of multiplicative type, and
therefore is central (\cite{milne2017}, 12.38). Thus $G^{\mathrm{ad}}$ is a
quotient of $(G_{\mathrm{Hg}}^{\mathrm{der}})^{\mathrm{ad}}=G_{\mathrm{Hg}%
}^{\mathrm{ad}}$, and so is a Mumford--Tate group. Conversely, suppose that
$G^{\mathrm{ad}}$ is a Mumford--Tate group, and let $\bar{h}\colon\mathbb{S}%
{}\rightarrow G_{\mathbb{R}}^{\mathrm{ad}}$ be a homomorphism satisfying
(MT1,2). According to Proposition \ref{m16}, there exists a reductive group
$G^{\prime}$ with $G^{\prime\mathrm{der}}=G$ and a homomorphism $h^{\prime
}\colon\mathbb{S}{}\rightarrow G_{\mathbb{R}{}}^{\prime}$ lifting $\bar{h}$
and satisfying (MT1,2). Let $\rho\colon G_{\mathrm{Hg}}\rightarrow G^{\prime}$
be the homomorphism given by Corollary \ref{m20}. Then $\rho$ maps
$(G_{\mathrm{Hg}})^{\mathrm{der}}$ onto $(G^{\prime})^{\mathrm{der}}=G$.

(c) Apply \ref{m10a} and \ref{m10b}.
\end{proof}

\begin{corollary}
\label{m23}The pro-algebraic group $G^{\mathrm{der}}$ is simply connected.
\end{corollary}

\begin{proof}
Immediate consequence of (b) of the theorem.
\end{proof}

\begin{remark}
\label{m19}The group $G_{\mathrm{Hg}}$ is not a product of $G_{\mathrm{Hg}%
}^{\mathrm{der}}$ and $S$ because this would imply that, for a semisimple
algebraic group $H$ over $\mathbb{Q}{}$, every homomorphism $h\colon
\mathbb{S}{}\rightarrow H_{\mathbb{R}}^{\mathrm{ad}}$ satisfying (MT1,2) lifts
to a homomorphism $\mathbb{S}{}\rightarrow H_{\mathbb{R}}$ satisfying the same
conditions, but this is not true.
\end{remark}

\begin{remark}
\label{m25}No description, even conjectural, is known for the essential image
of the functor from the category of motives over $\mathbb{C}$ to the category
of polarizable rational Hodge structures. However, when one replaces the
category of all motives with the subcategory generated by the motives of
abelian varieties such a description is known (Theorem 1.27 of
\cite{milne1994}). Note that these Hodge structures need be neither of CM-type
nor of weight 1. The proof of Theorem \ref{m17} is a (simpler) variant of the
proof of that theorem.
\end{remark}

\bibliographystyle{cbe}
\bibliography{CMT}

\end{document}